\documentclass{article}
\usepackage{authblk}

\usepackage{tikz}
\usepackage{booktabs}
\usepackage{algorithm}  
\usepackage[end]{algpseudocode}

\usepackage{graphicx}
\usepackage{xcolor}
\usepackage{epsfig} 
\usepackage{mathptmx} 
\usepackage{times} 
\usepackage{amsmath} 
\usepackage{amssymb}  
\usepackage[hidelinks]{hyperref}
\usepackage{lipsum}
\usepackage{mathtools}
\usepackage{cuted}
\usepackage{multirow}   
\usepackage{float}
\usepackage{blindtext}
\usepackage [list=true]{subcaption}
\usepackage{color}
\usepackage[numbers, sort & compress]{natbib}
\usepackage[doublespacing]{setspace}
\DeclareMathOperator{\Tr}{Tr}

\newtheorem{lemma}{Lemma}
\newtheorem{theorem}{Theorem}

\begin{document}

\title{Stochastic Programming with Primal-Dual Dynamics: A Mean-Field Game Approach}



\author[1]{C. T. R\"oling}
\author[2]{D. Bauso}
\author[3]{H. Tembine}

\affil[1]{ENTEG, Faculty of Science and Engineering, University of Groningen,
               Nijenborgh 4 9747 AG Groningen, The Netherlands,
               casperthomas.roling@gmail.com}
\affil[2]{Jan C. Willems Center for Systems and Control,
               ENTEG, Faculty of Science and Engineering,
                University of Groningen,
                 Nijenborgh 4 9747 AG Groningen, The Netherlands, and
              Dipartimento di Ingegneria,
              Universit\`a di Palermo,
              Palermo, Italy,
              d.bauso@rug.nl}

\affil[3]{
              Learning and Game Theory Laboratory,
               New York University Abu Dhabi,
               tembine@nyu.edu}


\maketitle

\begin{abstract}
This study addresses primal-dual dynamics for a stochastic programming problem for capacity network design. It is proven that consensus can be achieved on the \textit{here and now} variables which represent the capacity of the network. The main contribution is  a heuristic approach which involves the formulation of the problem as a mean-field game. Every agent in the mean-field game has control over its own primal-dual dynamics and seeks consensus with neighboring agents according to a communication topology. We obtain theoretical results concerning the existence of a mean-field equilibrium. Moreover, we prove that the consensus dynamics converge such that the agents agree on the capacity of their respective micro-networks. Lastly, we emphasize how penalties on control and state influence the dynamics of agents in the mean-field game.
\end{abstract}




\section{Introduction}
This work focuses on primal-dual dynamics for a stochastic programming problem related to network capacity design. Demand is extracted from sink nodes and goods flow over the edges. Edges are subject to capacity constraints and flow conservation constraints, implying that assigned capacity cannot be exceeded and that the network should comply to demand being extracted from sink nodes, respectively \cite{stolyar2006greedy, bertsimas2000traffic}. 
The setup fits a broad range of application cases such as communications network optimization, supply chain network optimization, energy resource allocation, (cyber-physical) smart grids, and stock-market pricing \cite{bauso2017consensus, feijer2010stability, chen2011convergence, ferragut2014network, li2017distributed,zhao2012swing, li2017dynamic, notarstefano2019distributed}. 
The call for more efficient supply chain logistics with an ever-growing culture of online ordering in addition to the upcoming era of the \textit{physical internet} explains the relevance for this study \cite{montreuil2011toward, montreuil2010towards}. Also, the current issues with climate change and Earth resource depletion again stress the need for efficient logistics \cite{halldorsson2010sustainable,piecyk2015green}. 


Recent works have been dedicated to stability analysis of primal-dual dynamics with non-strong convexity issues where is proven that under specific conditions, the primal-dual dynamics are globally exponentially stable  \cite{chen2019exponential, liang2019exponential}. Furthermore, analysis on robustness has been performed in \cite{nguyen2018contraction, cherukuri2017role, jonsson2010primal}, where the robustness of stability of the primal-dual dynamics is analyzed and performance guarantees are provided for different primal-dual systems, such as power systems. Also, the global/local nature of the convexity-concavity properties of the primal-dual dynamics was examined. Explicit attention has been paid to asymptotic stability properties of the primal-dual dynamics in \cite{feijer2010stability, cherukuri2016asymptotic, cherukuri2017saddle, stolyar2006greedy}. In these works, the well-known phenomena of Lyapunov functions and the invariance principle by LaSalle are used to conduct stability analyses. 

There is less work dedicated to the combination of primal-dual dynamics with mean-field games. However, there have been studies concerning mean-field games applied on a primal-dual problem, such as in \cite{achdou2012mean}, where a planning problem is discussed. Additionally, in \cite{yang2018distributed} and \cite{yang2016distributed} Lagrangian relaxation is used in mean-field game power control.

The main contribution is a heuristic approach which involves the reformulation of the original problem into a large number of stochastic primal-dual dynamics which are coupled in the same spirit as in mean-field games. To be more specific, we consider the instance of a large number of players, each assigned with a primal-dual dynamics subject to a specific realization of the demand. The players of the mean-field game have to obtain consensus over the \textit{here and now} decision variables: the capacity of each edge of the supply chain network. We also consider \textit{wait and see} decision variables for flow and the Lagrange multipliers coming from the primal-dual dynamics. 
 Note the analogy with a two-player zero-sum game whereby the first player, the \textit{minimizer}, sets the \textit{here and now} variables required to accommodate for any potential scenario concerning flow and demand realization. Subsequently, the second player, the \textit{maximizer}, sets the \textit{wait and see} variables to maximize demand satisfaction given the decision of player 1 and the realized demand. The condition for convergence to a consensus value is that the communication network of agents needs to be connected, i.e., for any pair of nodes there exists a path connecting them. We first transform the primal-dual dynamics into a mean-field game, after which theoretical results are provided concerning a mean-field equilibrium solution in the same spirit as in \cite{bauso2017consensus,bauso2016crowd,bauso2017dynamic}. Lastly, we provide numerical simulations to prove once again that consensus is obtained.  Since the proposed methodology is heuristic in nature, the convergence values are in general sub-optimal.

We wish to stress that  though the obtained mean-field game has a classic structure, the game is obtained via an original
methodology. The methodology involves turning the micro-network model into a set of primal-dual dynamics and after that into a mean-field game in which each agent faces a different realization of the uncertain demand $\omega$. We see this as a
value and believe that the link between primal-dual dynamics and mean-field games expands the significance,
meaning, and potential of mean-field games approaches beyond the ones already in the literature. Note also that the obtained mean-field game depends strongly on the micro-network optimization parameters. Actually the incidence matrix, the penalty coefficients and the uncertain demand enter into the mean-field dynamics of the optimization flow and capacity variables as well as of the Lagrange multipliers.

This paper is organized as follows. In Section~\ref{mf}, we formulate the stochastic programming for network design. In Section~\ref{pd}, we present the primal-dual dynamics. In Section~\ref{mfg}, present the corresponding mean-field game. In Section \ref{res}, we state theoretical results on mean-field equilibrium and convergence. In Section \ref{num}, we provide a simulation example to corroborate our results. Lastly, in Section \ref{conc}, we provide concluding remarks and discuss future works.

\section{Stochastic programming for network design}\label{mf}
\begin{figure} [h]
\centering
\begin{tikzpicture}
\begin{scope}[shift={(0,0)},scale=1]
\draw [black,fill=gray!20](0,-0.15) ellipse (0.18cm and .15cm);
\draw [black,fill=gray!20] (0,0) circle (.1cm); 
\draw (0.1,-0.15) .. controls (.12,-.16) and (0.14,-.18) .. (0.13,-0.25);
\draw (-0.1,-0.15) .. controls (-.12,-.16) and (-0.14,-.18) .. (-0.13,-0.25);
\end{scope}
\begin{scope}[shift={(1,-0.3)},scale=1]
\draw (0,-0.15) ellipse (0.18cm and .15cm);
\draw [fill=white] (0,0) circle (.1cm); 
\draw (0.1,-0.15) .. controls (.12,-.16) and (0.14,-.18) .. (0.13,-0.25);
\draw (-0.1,-0.15) .. controls (-.12,-.16) and (-0.14,-.18) .. (-0.13,-0.25);
\end{scope}
\begin{scope}[shift={(2,0)},scale=1]
\draw [black,fill=gray!20] (0,-0.15) ellipse (0.18cm and .15cm);
\draw [black,fill=gray!20] (0,0) circle (.1cm); 
\draw (0.1,-0.15) .. controls (.12,-.16) and (0.14,-.18) .. (0.13,-0.25);
\draw (-0.1,-0.15) .. controls (-.12,-.16) and (-0.14,-.18) .. (-0.13,-0.25);
\end{scope}
\begin{scope}[shift={(2.5,.5)},scale=1]
\draw [black,fill=gray!20] (0,-0.15) ellipse (0.18cm and .15cm);
\draw [black,fill=gray!20] (0,0) circle (.1cm); 
\draw (0.1,-0.15) .. controls (.12,-.16) and (0.14,-.18) .. (0.13,-0.25);
\draw (-0.1,-0.15) .. controls (-.12,-.16) and (-0.14,-.18) .. (-0.13,-0.25);
\end{scope}
\begin{scope}[shift={(2,1)},scale=1]
\draw (0,-0.15) ellipse (0.18cm and .15cm);
\draw [fill=white] (0,0) circle (.1cm); 
\draw (0.1,-0.15) .. controls (.12,-.16) and (0.14,-.18) .. (0.13,-0.25);
\draw (-0.1,-0.15) .. controls (-.12,-.16) and (-0.14,-.18) .. (-0.13,-0.25);
\end{scope}
\begin{scope}[shift={(0,1)},scale=1]
\draw (0,-0.15) ellipse (0.18cm and .15cm);
\draw [fill=white] (0,0) circle (.1cm); 
\draw (0.1,-0.15) .. controls (.12,-.16) and (0.14,-.18) .. (0.13,-0.25);
\draw (-0.1,-0.15) .. controls (-.12,-.16) and (-0.14,-.18) .. (-0.13,-0.25);
\end{scope}
\begin{scope}[shift={(1,1.3)},scale=1]
\draw (0,-0.15) ellipse (0.18cm and .15cm);
\draw [fill=white] (0,0) circle (.1cm); 
\draw (0.1,-0.15) .. controls (.12,-.16) and (0.14,-.18) .. (0.13,-0.25);
\draw (-0.1,-0.15) .. controls (-.12,-.16) and (-0.14,-.18) .. (-0.13,-0.25);
\end{scope}
\begin{scope}[shift={(-0.5,0.5)},scale=1]
\draw (0,-0.15) ellipse (0.18cm and .15cm);
\draw [fill=white] (0,0) circle (.1cm); 
\draw (0.1,-0.15) .. controls (.12,-.16) and (0.14,-.18) .. (0.13,-0.25);
\draw (-0.1,-0.15) .. controls (-.12,-.16) and (-0.14,-.18) .. (-0.13,-0.25);
\end{scope}
\begin{scope}[shift={(0.5,0.5)},scale=1]
\draw (0,-0.15) ellipse (0.18cm and .15cm);
\draw [fill=white] (0,0) circle (.1cm); 
\draw (0.1,-0.15) .. controls (.12,-.16) and (0.14,-.18) .. (0.13,-0.25);
\draw (-0.1,-0.15) .. controls (-.12,-.16) and (-0.14,-.18) .. (-0.13,-0.25);
\end{scope}
\begin{scope}[shift={(1.5,0.5)},scale=1]
\draw (0,-0.15) ellipse (0.18cm and .15cm);
\draw [fill=white] (0,0) circle (.1cm); 
\draw (0.1,-0.15) .. controls (.12,-.16) and (0.14,-.18) .. (0.13,-0.25);
\draw (-0.1,-0.15) .. controls (-.12,-.16) and (-0.14,-.18) .. (-0.13,-0.25);
\end{scope}
\draw [very thick, gray] (0.2,-0.2)--(0.8,-0.4);\draw [very thick, gray] (1.2,-0.4)--(1.8,-0.2);
\draw [very thick, gray] (2.15,-0.05)--(2.4,0.2);\draw [very thick, gray] (2.4,0.58)--(2.2,0.75);
\draw [very thick, gray] (1.85,1.00)--(1.2,1.15);\draw [very thick, gray] (0.8,1.15)--(0.2,1.00);
\draw [very thick, gray] (-0.2,0.80)--(-0.4,0.60);\draw [very thick, gray] (-0.4,0.20)--(-0.2,0.00);
\draw [very thick, gray] (0.1,0.10)--(0.3,0.3);\draw [very thick, gray] (1.9,0.10)--(1.7,0.25);
\draw [very thick, gray] (0.7,0.40)--(1.3,0.4);
\draw [very thick, gray] (1.6,0.60)--(1.8,0.8);
\draw [very thick, gray] (0.4,0.60)--(0.2,0.8);
\begin{scope}[shift={(-1.0,-2.5)},scale=1]
\draw [fill=white] (0,0) circle (.1cm); \draw [fill=white] (1,0) circle (.1cm); 
\draw [fill=white] (-0.5,0.5) circle (.1cm); \draw [fill=white] (1.5,0.5) circle (.1cm); 
\draw [fill=white] (-0.5,-0.5) circle (.1cm); \draw [fill=white] (1.5,-0.5) circle (.1cm); 
\draw [->] (-1,-0.5)--(-0.6,-0.5);\draw [->] (-1,0.5)--(-0.6,0.5);
\draw [->] (-0.4,-0.5)--(1.4,-0.5);\draw [->] (-0.4,0.5)--(1.4,0.5);
\draw [->] (1.6,-0.5)--(2.0,-0.5);\draw [->] (1.6,0.5)--(2.0,0.5);
\draw [->] (-0.4,-0.4)--(-0.1,-0.1);\draw [->] (-0.4,0.4)--(-0.1,0.1);
\draw [->] (0.1,0.0)--(0.9,0.0);
\draw [->] (1.1,0.1)--(1.4,0.4);\draw [->] (1.1,-0.1)--(1.4,-0.4);
\draw [fill=gray,opacity = 0.3] (0.5,0) ellipse (2.0cm and 1.2cm);
\draw  (-1,.9)--(0.8,2.2);\draw  (2.1,.8)--(1.2,2.1);
\end{scope}
\begin{scope}[shift={(1.9,-1.45)},scale=.5]
\draw [fill=white] (0,0) circle (.1cm); \draw [fill=white] (1,0) circle (.1cm); 
\draw [fill=white] (-0.5,0.5) circle (.1cm); \draw [fill=white] (1.5,0.5) circle (.1cm); 
\draw [fill=white] (-0.5,-0.5) circle (.1cm); \draw [fill=white] (1.5,-0.5) circle (.1cm); 
\draw [->] (-1,-0.5)--(-0.6,-0.5);\draw [->] (-1,0.5)--(-0.6,0.5);
\draw [->] (-0.4,-0.5)--(1.4,-0.5);\draw [->] (-0.4,0.5)--(1.4,0.5);
\draw [->] (1.6,-0.5)--(2.0,-0.5);\draw [->] (1.6,0.5)--(2.0,0.5);
\draw [->] (-0.4,-0.4)--(-0.1,-0.1);\draw [->] (-0.4,0.4)--(-0.1,0.1);
\draw [->] (0.1,0.0)--(0.9,0.0);
\draw [->] (1.1,0.1)--(1.4,0.4);\draw [->] (1.1,-0.1)--(1.4,-0.4);
\draw [fill=gray,opacity = 0.3] (0.5,0) ellipse (2.0cm and 1.2cm);
\draw  (-1,.9)--(0.0,2.2);\draw  (2.1,.8)--(0.6,2.1);
\end{scope}
\begin{scope}[shift={(3.5,-0.5)},scale=.4]
\draw [fill=white] (0,0) circle (.1cm); \draw [fill=white] (1,0) circle (.1cm); 
\draw [fill=white] (-0.5,0.5) circle (.1cm); \draw [fill=white] (1.5,0.5) circle (.1cm); 
\draw [fill=white] (-0.5,-0.5) circle (.1cm); \draw [fill=white] (1.5,-0.5) circle (.1cm); 
\draw [->] (-1,-0.5)--(-0.6,-0.5);\draw [->] (-1,0.5)--(-0.6,0.5);
\draw [->] (-0.4,-0.5)--(1.4,-0.5);\draw [->] (-0.4,0.5)--(1.4,0.5);
\draw [->] (1.6,-0.5)--(2.0,-0.5);\draw [->] (1.6,0.5)--(2.0,0.5);
\draw [->] (-0.4,-0.4)--(-0.1,-0.1);\draw [->] (-0.4,0.4)--(-0.1,0.1);
\draw [->] (0.1,0.0)--(0.9,0.0);
\draw [->] (1.1,0.1)--(1.4,0.4);\draw [->] (1.1,-0.1)--(1.4,-0.4);
\draw [fill=gray,opacity = 0.3] (0.5,0) ellipse (2.0cm and 1.2cm);
\draw  (-1.5,.8)--(-1.9,1.9);\draw  (1.4,1.2)--(-1.6,1.9);
\end{scope}
\end{tikzpicture}
\caption{Two-layer network: macro-network $\hat G(\hat V,\hat E)$ in layer~1 and micro-network $G(V,E)$ in layer~2.}  \label{fig:twolayernetwork}
\end{figure}
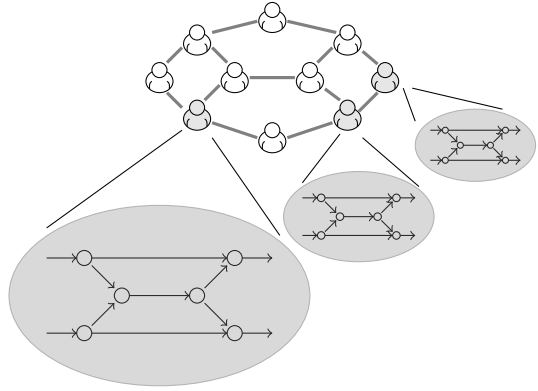

\begin{figure} [h]
\centering
\begin{tikzpicture}
\begin{scope}[shift={(-1.0,-2.5)},scale=2.5]
\draw [very thick] (0,0) circle (.1cm) ;\node at (0,0){$5$};
 \draw [very thick] (1,0) circle (.1cm); \node at (1,0){$6$};
\draw [very thick] (-0.5,0.5) circle (.1cm);\node at (-0.5,0.5){$1$};
 \draw [very thick] (1.5,0.5) circle (.1cm); \node at (1.5,0.5){$4$};
\draw [very thick]  (-0.5,-0.5) circle (.1cm); \node at (-0.5,-0.5){$2$};
\draw [very thick] (1.5,-0.5) circle (.1cm); \node at (1.5,-0.5){$3$};
\draw [->,very thick] (-1,-0.5)--(-0.6,-0.5);\node at (-0.8,-0.4){$u_2$};\node at (-0.8,-0.6){$c_2$};
\draw [->, very thick] (-1,0.5)--(-0.6,0.5);\node at (-0.8,0.6){$u_1$};\node at (-0.8,0.4){$c_1$};
\draw [->,very thick] (-0.4,-0.5)--(1.4,-0.5);\node at (0.5,-0.4){$u_6$};\node at (0.5,-0.6){$c_6$};
\draw [->, very thick] (-0.4,0.5)--(1.4,0.5);\node at (0.5,0.6){$u_3$};\node at (0.5,0.4){$c_3$};
\draw [->, very thick] (1.6,-0.5)--(2.0,-0.5);\node at (1.8,-0.4){$\omega_3$};
\draw [->,very thick] (1.6,0.5)--(2.0,0.5);\node at (1.8,0.6){$\omega_4$};
\draw [->,very thick] (-0.4,-0.4)--(-0.1,-0.1);\node at (-0.1,-0.3){$c_5$};\node at (-0.3,-0.1){$u_5$};
\draw [->,very thick] (-0.4,0.4)--(-0.1,0.1);\node at (-0.3,0.1){$u_4$};\node at (-0.1,0.3){$c_4$};
\draw [->,very thick] (0.1,0.0)--(0.9,0.0);\node at (0.5,0.1){$u_8$};\node at (0.5,-0.1){$c_8$};
\draw [->,very thick] (1.1,0.1)--(1.4,0.4);\node at (1.1,0.3){$u_9$};\node at (1.3,0.1){$c_9$};
\draw [->,very thick] (1.1,-0.1)--(1.4,-0.4);\node at (1.1,-0.3){$u_7$};\node at (1.3,-0.1){$c_7$};
\end{scope}
\end{tikzpicture}
 \caption{A graphical representation of the directed network corresponding to the first problem instance.}
    \label{fig:firstnetwork}
\end{figure}
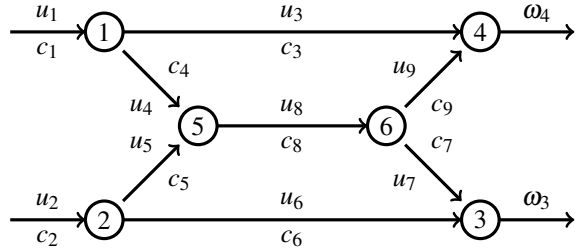

Consider a two-layer network as in Fig.~\ref{fig:twolayernetwork}. Layer 1 involves the communication topology (henceforth referred to as \emph{macro-network}) denoted as $\hat G(\hat V,\hat E)$ whereas layer 2 comprises the physical flow network (henceforth referred to as \emph{micro-network}) denoted as $G(V,E)$. We first define the micro-network as follows. Let the set of nodes be $V \in \{V_1, V_2,\ldots, V_n\}$. Parameter $n$ denotes the number of nodes present in the network. Let  the set of edges be $E \in \{E_1, E_2,\ldots,E_m\}$. Parameter $m$ denotes the number of edges in the network. Demand $\omega$ pulling from sink nodes is described by a vector $\omega \in \mathbb{R}^{n\times1}$. Capacity of the edges are the \textit{here and now} variables and are denoted by $c \in \mathbb{R}^{m\times1}$. The transported goods in each edge are the \textit{wait and see} variables and are denoted by $u \in \mathbb{R}^{m\times1}$ denotes . A graphical representation of the micro-network can be found in Fig.~\ref{fig:firstnetwork}.

Now consider that all demand should always be fulfilled and backlog is discouraged in the system. Logically, the amount of transported goods cannot exceed the capacity for respective edges. 
The deterministic version of the problem we wish to solve involves the minimization over the flow $u$ and capacity $c$ of the edges and can be formulated as:
\begin{align}
\begin{split} \label{eq:fvalstochastic}
\min \quad f_1(c)+f_2(u) \\
\begin{aligned} 
\textup{s.t.} \quad  \quad \tilde Bu &=\omega, \\
                   \quad u  &\leq  c, \\
         u,c  &>  0, 
 \end{aligned}
\end{split}
\end{align} 
where $\tilde B\in \mathbb{R}^{n\times m}$ is the incidence matrix of the micro-network. In the optimization problem, we let $f_1(c)=\tfrac{1}{2}c^T\widetilde{Q}_1c+\widetilde{f}_1^Tc$, and analogously, $f_2(u)=\tfrac{1}{2}u^T\widetilde{Q}_2u+\widetilde{f}_2^Tu$. Lastly, let us define $\widetilde{Q}_1,\widetilde{Q}_2 \in \mathbb{R}^{m\times m}$.

Under the assumption that $\omega$ is an uncertain random parameter, we may model this problem as a two-stage stochastic program as follows:
\begin{align}
    \begin{split}
\min \quad f_1(c)+\mathbb{E}_\omega\hat{Q}(u,\omega) \\
\begin{aligned}
\textup{s.t.} \quad  \quad \tilde Bu &=\omega,\\
                  \quad u  &\leq  c, \\
         u,c  &>  0. 
\end{aligned} 
    \end{split}
\end{align}
Here, $\mathbb E_\omega(.)$ denotes expectation with respect to $\omega$ and $\hat{Q}(u,\omega)$ denotes the sub-optimal value of the second-stage problem, namely  $$\hat{Q}(u,\omega):=\min_u\{f_2(u)|Bu =\omega \}.$$

\section{Primal-dual dynamics}\label{pd}

From the Lagrangian and, after rearranging, we obtain:
\begin{equation} \nonumber
    \begin{split}
        &\mathcal{L}(u,c,\lambda,\mu)  \\
        &= f_1(c)+\lambda^T(\tilde Bu-\omega)+\mu^T(u-c) +\mathbb{E}f_2(u) \\
        & = f_1(c)+\left[\tilde B^T\lambda+\mu \right]^Tu-\mu^Tc-\lambda^T\omega+\mathbb{E}f_2(u) \\
        & = f_1(c) + \left[\begin{bmatrix}I\\ -I \end{bmatrix}  \mu + \begin{bmatrix} \tilde B^T \\ 0 \end{bmatrix}\lambda \right]^T \begin{bmatrix}u\\c \end{bmatrix}-\omega^T\lambda+\mathbb{E}f_2(u), 
    \end{split}
\end{equation} 
where  $\mu \in \mathbb{R}^{m\times1}$ and $\lambda \in \mathbb{R}^{n\times1}$ are the Lagrange multipliers. 

 The Lagrange dual function can then be specified as follows:
\begin{equation}\nonumber
\begin{split}
    \inf_{u,c}\Biggl\{ f_1(c) + \left[\begin{bmatrix}I\\ -I \end{bmatrix}  \mu + \begin{bmatrix} \tilde B^T \\ 0 \end{bmatrix}\lambda \right]^T \begin{bmatrix}u\\c \end{bmatrix} 
     -\omega^T\lambda+\mathbb{E}f_2(u) \Biggr\}
    \end{split}.
\end{equation} 
For a solution in nonlinear programming to be optimal, the Karush-Kuhn-Tucker (KKT) conditions must hold, which we specify next: 
\begin{equation}
    \begin{bmatrix} \nabla_uf_2(u) \\ \nabla_cf_1(c) \end{bmatrix}+ \begin{bmatrix} I\\-I \end{bmatrix}\mu+\begin{bmatrix} \tilde B^T\\0\end{bmatrix}\lambda = 0.
\end{equation} 
The gradient of the Lagrangian with respect to $u$ and $c$ vanishes, which therefore results in the following KKT conditions:
\begin{equation}
    \begin{cases}
    \mu+\tilde B^T\lambda+\nabla_uf_2(u)=0, \\
    \nabla_cf_1(c)-\mu =0.
    \end{cases}
\end{equation} 

We formulate the primal-dual dynamics according to gradient descent and gradient ascent for convex and concave variables, respectively:
\begin{equation}
    \begin{cases}
    \begin{aligned}
    \dot{u}&=\Bigl[-\nabla_u\mathcal{L}\Bigr]_+=\Bigl[\nabla_uf_2(u)-I\mu- \tilde B^T\lambda\Bigr]_+\\
    &=\Bigl[\mu- \tilde B^T\lambda-\widetilde{Q}_2u-\widetilde{f}_2\Bigr]_+ \\
    \dot{c}&= \Bigl[-\nabla_c\mathcal{L}\Bigr]_+=\Bigl[\nabla_uf_1(c)-(-I)\mu\Bigr]_+\\
    & =\Bigl[-\widetilde{Q}_1c-\widetilde{f}_1+\mu \Bigr]_+ \\
    \dot{\lambda}&= +\nabla_\lambda\mathcal{L}=\tilde Bu-\omega \\
   \dot{\mu}&= \Bigl[+\nabla_\mu\mathcal{L}u-c\Bigr]_+
    \end{aligned}
    \end{cases}
\end{equation} 
In compact vector form, the system can be written as 
\begin{equation}
\begin{split}
    \begin{bmatrix} \dot{u}\\\dot{c}\\\dot{\lambda}\\\dot{\mu}\end{bmatrix}&=\begin{bmatrix} -\widetilde{Q}_2&0&-\tilde B^T&-I\\ 0&-\widetilde{Q}_1&0&I\\B&0&0&0\\I&-I&0&0 \end{bmatrix}\begin{bmatrix}u\\c\\\lambda\\\mu \end{bmatrix}+\begin{bmatrix}-\widetilde{f}_2\\-\widetilde{f}_1\\-\omega\\0 \end{bmatrix} \\
    & \in \mathbb{R}^m_+ \times \mathbb{R}^m_+ \times \mathbb{R}^n \times \mathbb{R}^m_+.
    \label{eq:modelcompactform}
    \end{split}
\end{equation} 
As previously explained, we would like to obtain consensus regarding the capacity of edges of the connected undirected graph $G=(V,E)$ defined earlier. This is possible formulating the problem as a mean-field game as explained next. 


\section{Construction of the Mean-Field Game}\label{mfg}

Consider $p$ populations such that each player of our game belongs to a population $k\in \left \{ 1,\ldots,p\right \}$. The generic player of population $k$ is characterized by its state $x(t) \in \mathbb{R}^{3m+n}$ at time $t$ for a given time horizon window $[0,T]$ which evolves according to (\ref{eq:modelcompactform}) under the control variable $v_k$. In particular,  
\begin{equation}
\begin{split}
    \underbrace{\begin{bmatrix} \dot{u}\\\dot{c}\\\dot{\lambda}\\\dot{\mu}\end{bmatrix}}_{\dot{x}}&=\underbrace{\begin{bmatrix} -\widetilde{Q}_2&0&-\tilde B^T&-I\\ 0&-\widetilde{Q}_1&0&I\\\tilde B&0&0&0\\I&-I&0&0 \end{bmatrix}}_A\underbrace{\begin{bmatrix}u\\c\\\lambda\\\mu \end{bmatrix}}_x+\underbrace{\begin{bmatrix}-\widetilde{f}_2\\-\widetilde{f}_1\\-\omega\\0 \end{bmatrix}}_C \\
    & + \underbrace{\begin{bmatrix} 0 & \dots & \dots & 0 \\ \vdots & \mathbf{1} & \dots & \vdots\\ \vdots & \vdots & \ddots &\vdots \\ 0 & \dots & \dots & 0 \end{bmatrix}}_Bv_k.
    \end{split}
    \label{eq:modelcompactform2}
\end{equation}
Note that \eqref{eq:modelcompactform2} can be written in compact form as  $\dot x=Ax+Bv_k+C=:f(x,v_k)$.  Note that $x$ is the state vector of the system that involves also the decision variables of our original stochastic programming problem.


Now let us consider a probability density function that describes the density of the players of a population in state $x$ at time $t$, $m_k(x,t)$, with the property $\int_{\mathbb{R}^{3m+n}} m_k(x,t)dx=1$. 
Then the mean states can be computed following $\overline{m}_k(t)=\int_{\mathbb{R}^{3m+n}} xm_k(x,t)dx$.

To describe the interaction between populations, let us associate population $k$ with  agent $k$. Then we can introduce an interaction topology between agents, say $\hat G=\{\hat V,\hat E\}$, and we can define the neighbors of an agent $k$ in  $\hat G$ as:
$$
N(k) = \left \{ j \in \hat V \mid (k,j) \in \hat E \right \}.
$$
We define a player's objective in population $k$ based on the aggregate $k$th state as:
\begin{equation}
\rho_k=\frac{\sum_{j\in N(k)}\overline{m}_j}{|N(k)|}. \label{eq:rhok}
\end{equation}
Here, $|N(k)|$ denotes the cardinality of neighbor set $k$. Let us now consider a running cost function $g(x,\rho_k,v)$ and a terminal cost function $\psi(\rho_k,x)$, which are defined as follows:
\begin{equation*}
\begin{split}
g(x,\rho_k,v)&=\frac{1}{2}(\rho_k-x)^TQ(\rho_k-x)+\frac{1}{2}v^TRv,\\
\psi(\rho_k,x)&=\frac{1}{2}S(\rho_k-x)^2.
\end{split}
\end{equation*}
The first term in the expression for $g$ assigns a penalty to the state deviating from mean $\rho_k$. The second term assigns a penalty on control. $Q$, $S$ and $R$ are diagonal matrices of compatible dimensions. 

Every player in population $k$ wishes to solve the following problem:
\begin{equation}
    \min_{v(\cdot)}\mathbb{E}\int_0^T[g(x,\rho_k,v)]dt + \psi(\rho_k(T),X(T))
\end{equation}
subject to
$$
\dot{x}=Ax+Bv+C.
$$
For every population $k$, denote the value of the robust optimization problem starting at time $t$ and state $x$ by $\sigma_k(x,t)$. This results in the following mean-field game in $\sigma_k(x,t)$ and $m_k(x,t)$:
\begin{equation}
    \begin{cases}
    \partial_t\sigma_k(x,t)+\left \{ f(x,v^*)^T\partial_x\sigma_k(x,t)+g(x,\rho_k,v^*)\right \} = 0, \\
    \sigma_k(x,T) = \psi(\rho_k(T),x), \\
    \partial_tm_k(x,t)+div(m_k(x,t)f(\cdot))=0, \\
    m_k(x,0)=m_{k0}(x).
    \end{cases} 
    \label{eq:meanfieldgame}
\end{equation}
Any solution of \eqref{eq:meanfieldgame} is referred to as the \textit{mean-field equilibrium}, which provides the sub-optimal values for the \textit{wait and see} variables. Note that the second and fourth equation from \eqref{eq:meanfieldgame} are the boundary conditions, while the third equation is the advection equation. 

The optimal time-varying state-feedback control can be computed for every single agent in population $k$ and is given by:
\begin{equation}
\begin{split}
    v^*_k(x,t)& \in \arg \min_v \bigl\{ (Ax+Bv+C)^T\partial_x\sigma_k(x,t) \\
    &+g(x,\rho_k(t),v) \bigr\}.
    \end{split}
\end{equation}
In this expression, note that the Hamiltonian appears as the argument of the minimizer. 

\section{Mean-field equilibrium and convergence}\label{res}
In this section, we obtain an expression for the mean-field control and provide results for the mean-field equilibrium dynamics.

\begin{lemma}
The robust mean-field game takes the form:
\begin{equation}
    \begin{cases}
    \partial_t\sigma_k(x,t)-\frac{1}{2}(\partial_x\sigma_k(x,t))^T \left[BR^{-1}B^T \right]  \partial_x\sigma_k(x,t)\\ +(\partial_x\sigma_k(x,t))^T(Ax+C) \\
    +\frac{1}{2}(\rho_k-x)^TQ(\rho_k-x)=0 \\
    \sigma_k(x,T)=\psi(\rho_k(T),x) \\
    \partial_tm_k(x,t)\\+\partial_x\left[m_k(x,t)(Ax-BR^{-1}B^T\partial_x\sigma_k(x,t) +C) \right]=0 \\
    m_k(x,0)=m_{k0}(x) \\
    \overline{m}_k:=\int_\mathbb{R}xm_kdx
    \label{eq:systemth1}
    \end{cases}
\end{equation}
Additionally, the optimal control is:
\begin{equation}
    v_k^*=-R^{-1}B^T\partial_x\sigma_k(x,t). \label{eq:optcontrol}
\end{equation}
\end{lemma}
In this set of equations, the first equation corresponds to the Hamilton-Jacobi-Isaacs equation, and the third equation corresponds to the Fokker-Planck-Kolmogorov equation. The proof of \textit{Lemma 5.1} can be found in Appendix A.

We assume that the time evolution of the common state is known and, subsequently, investigate the solution of the Hamilton-Jacobi equation. Consider the following problem with known $\rho_k$:
\begin{equation}
    \min_{v(\cdot)}\mathbb{E}\int_0^T \left[ g(X(t),\rho_k(t),v(t)) \right]dt
\end{equation}
where
$$
\dot{x}=Ax+Bv+C.
$$
The next theoretical result presents the mean-field equilibrium control.
 In preparation to that, let us consider a probability density function that describes the density of the players of a population in state $c$ at time $t$, $m_k(c,t)$, with the property $\int_{\mathbb{R}^{m}} m_k(c,t)dc=1$. 
Then the mean states can be computed following $\overline{m}_k^c(t)=\int_{\mathbb{R}^{m}} c m_k(c,t)dc$.

We define a player's objective in population $k$ based on the aggregate $k$th state as:
\begin{equation}
\rho_k^c=\frac{\sum_{j\in N(k)}\overline{m}_j^c}{|N(k)|}. \label{eq:rhok}
\end{equation}

In the following, given a generic matrix $A  \in \mathbb{R}^{(3m+n) \times (3m+n)}$ we denote by $A_c  \in \mathbb{R}^{m}$ the matrix obtained extracting the rows and columns associated with the only variables $c \in \mathbb{R}^{m \times m}$. In addition,  
we consider the following value function:
\begin{equation}
    \sigma_k(x,t)=\frac{1}{2}x^T\Phi(t)x+H(t)^Tx+\chi(t).
\end{equation}
and denote 
$$
\widetilde{Z}:=\left[A^T-2\Phi BR^{-1}B^T \right]^{-1}.
$$
\begin{theorem}\label{thm1}
A mean-field equilibrium for the dynamics of \eqref{eq:systemth1} is obtained from the following set of equations:
\begin{equation}
\begin{cases}
    \sigma_k(x,t)=\frac{1}{2}x^T\Phi(t)x+H(t)^Tx+\chi(t), \\
    \begin{aligned}
    \dot{\overline{m}}_k(t)&=\left[A-BR^{-1}B^T\Phi(t)\right]\overline{m}_k(t)\\
    &-BR^{-1}B^TH(t)+C, 
    \end{aligned}
    \label{eq:th2mfe}
\end{cases}
\end{equation}
where
\begin{equation}
    \begin{cases}
    \begin{aligned}
    \dot{\Phi}(t)&+A^T\Phi(t)+\Phi^T\left[-BR^{-1}B^T\right]\Phi(t)+Q=0 \\
    &\in [0,T[,\quad \Phi(T)=S, \end{aligned}\\
    \begin{aligned}
    \dot{H}(t)&-2\Phi(t) BR^{-1}B^T H(t)+A^TH(t)\\
    &+\Phi(t)^TC-Q\rho_k(t)=0 \in [0,T[,\\
    &\quad h(T) = -S\rho_k(T), \end{aligned} \\
    \begin{aligned}\dot{\chi}&+H(t)^T\left[BR^{-1}B^T\right]H(t)+H(t)^TC\\
    &+\frac{1}{2}\rho_k(t)^TQ\rho_k(t)\in [0,T[,\quad \chi(T)=\frac{1}{2}\rho_k^T(T)S\rho_k(T).
    \end{aligned}
    \label{eq:th2phihchi}
    \end{cases}
\end{equation}
Additionally, the mean-field equilibrium control is:
\begin{equation}
    v^*_k=-R^{-1}B^T(\Phi(t)^Tx+H(t)).
\end{equation}

 Furthermore, for infinite time horizon $T \rightarrow \infty$, and for all aggregate states $\overline{m}^c = \left( \overline{m}_1^c,\overline{m}_2^c,\ldots,\overline{m}_p^c \right)$, we have the following consensus-type dynamics:
\begin{equation}
\dot{\overline{m}}^c = - \mathbf{\tilde Q_1} \overline{m}^c(t) -L \overline{m}^c(t) +\delta.
\end{equation}
where, $\mathbf{\tilde Q_1}:=diag(\tilde Q_1)$ (diagonal matrix with block entry $\tilde Q_1$), $\delta:=R_c^{-1} \Big( \widetilde{Z}_c \Phi_c^T \mu + (-\widetilde{Z}_cQ_c - \Phi_c^T) \rho_c^k \Big) 
+(\mu-\tilde f_1)$ and $L$ is defined as the graph-Laplacian matrix where the $kj$th entry is the block matrix:
\begin{equation}
L_{kj} =
    \begin{cases}
    R_c^{-1} \Phi_c^T  & \quad j=k, \\
\frac{1}{|N(k)|}R_c^{-1} \Phi_c^T & \quad j \neq k, j \in N(k), \\
    0 & \quad otherwise.
    \end{cases} 
 \end{equation}

%
%
\end{theorem}
This result is relevant since we can solve \eqref{eq:th2mfe} in closed form and we show that consensus is ultimately achieved among agents. The proof of Theorem \ref{thm1} can be found in Appendix B. As mentioned previously, the proposed methodology is heuristic in nature, and therefore the convergence values are in general sub-optimal.

\section{Simulation Example}\label{num}
\begin{figure}
    \centering
    \includegraphics[width=\columnwidth]{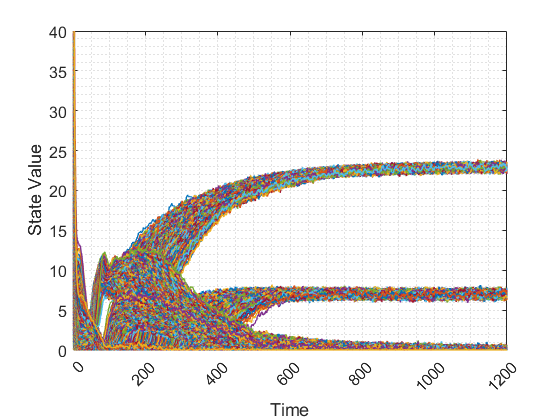}
    \caption{Simulation of the mean-field game with $Q=R=1$.}
    \label{fig:sim_QReq}
\end{figure}

\begin{figure}
    \centering
    \includegraphics[width=\columnwidth]{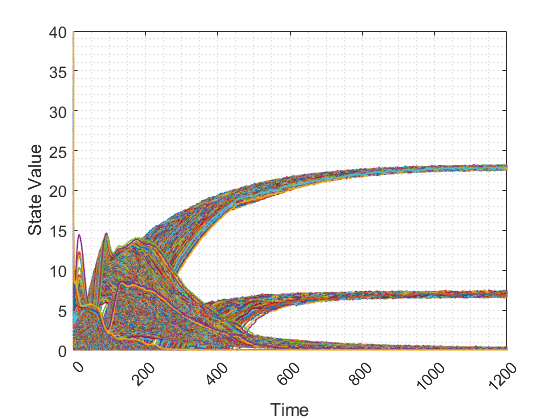}
    \caption{Simulation of the mean-field game with $Q=10$, $R=1$.}
    \label{fig:sim_Q10}
\end{figure}

\begin{figure}
    \centering
    \includegraphics[width=\columnwidth]{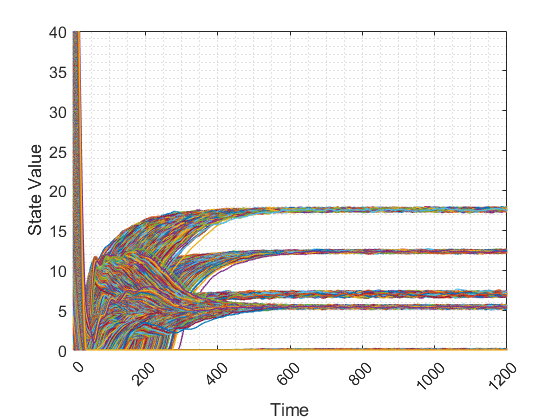}
    \caption{Simulation of the mean-field game with $Q=1$, $R=10$.}
    \label{fig:sim_R10}
\end{figure}

We are considering a scale-free network with $p=1000$, implying that the network consists of 1000 players. Players only consider the mean-field computed over its neighbors according to the expression of \eqref{eq:rhok}. 
Furthermore, we generate stochastic demand according to the normal distribution at every time instance for the sink nodes of the micro-networks. This makes the constant term vector $C$ change over time since $\omega$ is now renewed every time instance. Prior to the simulation, we parametrize the system as follows:
$$
\omega = \begin{bmatrix}0 & 0 & 23 & 7 & 0 & 0\end{bmatrix}^T,\quad \widetilde{Q}_1 = I^{9\times9},
$$
$$
\widetilde{Q}_2 = I^{9\times9},\quad \widetilde{f}_1 = \begin{bmatrix}1 & 1 & 1 & 1 & 1 & 1 & 1 & 1 & 1\end{bmatrix}^T,$$
$$
\widetilde{f}_2 = \begin{bmatrix}1 & 1 & 1 & 1 & 1 & 1 & 1 & 2 & 1 \end{bmatrix}^T,
$$
$$
C = \begin{bmatrix}-\widetilde{f}_2 & -\widetilde{f}_1 & -\omega & \mathbb{O}^{9\times1} \end{bmatrix}^T,\quad B=\begin{bmatrix} \mathbb{O}^{9\times1} \\ \mathbf{1}^{9\times1} \\ \mathbb{O}^{15\times1} \end{bmatrix},
$$
$$
B_I = \begin{bmatrix} 1 & 0 & -1 & -1 & 0 & 0 & 0 & 0 & 0 \\ 0 & 1 & 0 & 0 & -1 & -1 & 0 & 0 & 0 \\ 0 & 0 & 0 & 0 & 0 & 1 & 1 & 0 & 0 \\ 0 & 0 & 1 & 0 & 0 & 0 & 0 & 0 & 1 \\ 0 & 0 & 0 & 1 & 1 & 0 & 0 & -1 & 0 \\ 0 & 0 & 0 & 0 & 0 & 0 & -1 & 1 & -1
\end{bmatrix},$$
$$A = \begin{bmatrix}-\widetilde{Q}_2 & \mathbb{O}^{9\times9} & -B_I^T & -I^{9\times9} \\ \mathbb{O}^{9\times9} & -\widetilde{Q}_1 & \mathbb{O}^{9\times6} & I^{9\times9} \\ B_I & \mathbb{O}^{6\times9} & \mathbb{O}^{6\times6} & \mathbb{O}^{6\times9} \\ I^{9\times9} & -I^{9\times9} & \mathbb{O}^{9\times6} & \mathbb{O}^{9\times9} \end{bmatrix},
$$
$$
R = 1,\quad Q \in \mathbb{R}^{33\times33} = \begin{bmatrix} \mathbb{O}^{9\times9} & \dots & \dots & \mathbb{O}^{9\times9} \\ \vdots & \mathbf{1}^{9\times9} & \dots & \vdots\\ \vdots & \vdots & \ddots &\vdots \\ \mathbb{O}^{9\times9} & \dots & \dots & \mathbb{O}^{9\times9} \end{bmatrix}.
$$
 Note that in our example, the penalty of using flow $u_8$ is doubled since we want to discourage coupling. In general, any route can be incentivised by adjusting penalties and this emphasizes the versatility of our approach. Furthermore, in Table \ref{tab:simpar} the following additional simulation parameters are defined: Time step $\delta t$, mean and standard deviation of initial states $\widetilde{\sigma}_0$ and $\widetilde{\mu}_0$, and lastly, mean and standard deviation of demand at sink nodes 3 and 4 $\widetilde{\sigma}_{3,4}$ and  $\widetilde{\mu}_{3,4}$. The multi-population considered allows the simulation algorithm to be initialized with different initial conditions and also allows to accommodate different realizations of the \textit{wait and see} variables.
\begin{table}
   \caption{Simulation parameters and their respective values.} 
   \small 
   \centering 
   \begin{tabular}{lcccccccc} 
   \toprule[\heavyrulewidth]\toprule[\heavyrulewidth]
   \textbf{Parameter} & $\delta t$ & $\widetilde{\sigma}_{0}$ & $\widetilde{\mu}_0$ & $\widetilde{\mu}_3$ & $\widetilde{\mu}_4$ & $\widetilde{\sigma}_3$ & $\widetilde{\sigma}_4$ \\ 
   \midrule
   \textbf{Value} & 0.1 & 15 & 40 & 23 & 7 & 1 & 1 \\
   \bottomrule[\heavyrulewidth]
   \end{tabular} \label{tab:simpar}
\end{table}
We use the parameters determined as input for the following simulation algorithm by which we compute the states of each agent over time.
\begin{algorithm}
\caption{Simulation algorithm} \label{alg:simalg}
\begin{algorithmic}[1]
\State Define simulation parameters according to Table \ref{tab:simpar}
\State Define input matrices and vectors $\omega$, $\widetilde{f}_1$, $\widetilde{f}_2$, $\widetilde{Q}_1$, $\widetilde{Q}_2$, $B_I$, $A$, $Q$, $B$ and scalar $R$.
\State Compute initial state $x_0$
\State Initialize matrices $\rho$, $v^*$, $\omega$, and $C$
\State Solve the continuous time Riccati equation with $A,B,Q,R$ for $\Phi$
\State Compute $A_{ls}$ matrix for solving the system of linear equations for $H$
\For{Every time instance} 
    \For{Every agent}
        \State Compute average state of neighbors $\rho$
        \State Compute demand $\omega$ by the random normal distribution
        \State Compute vector $C$ with new vector $\omega$
        \State Compute $B_{ls}$ vector for solving the system of linear equations for $H$ 
        \State Solve for $H$
        \State Compute optimal control $v^*$
        \State Compute new state $x_{t+1}$
    \EndFor
\EndFor    
\end{algorithmic}
\end{algorithm}
In Figures \ref{fig:sim_QReq}, \ref{fig:sim_Q10}, and \ref{fig:sim_R10} is depicted how the states of the capacities evolve over time for all agents. Figure \ref{fig:sim_QReq} shows a simulation when the penalty on control equals the penalty on state deviation, i.e., $Q=R=1$. Figure \ref{fig:sim_Q10} presents a simulation of the mean-field game where the penalty on state deviation is increased, implying $Q=10$, $R=1$. Contrarily, the results of a mean-field game simulation with $Q=1$, $R=10$ can be found in Figure \ref{fig:sim_R10}.

By setting the penalty on state deviation sufficiently high as was done in the simulation of Figure \ref{fig:sim_Q10}, the converged states approximate each other to a larger extent compared to the case of equal $Q$ and $R$ because in this case the system is penalized more if states deviate from objective $\rho_k$. \\
Conversely, by setting the penalty on control sufficiently high, convergence is obtained quicker in comparison to the simulation with equal $Q$ and $R$. In this situation, the system is penalized more to make adjustments, i.e. control, hence the \textit{thinner} lines in the converged region. Another observation is that, in this case, the best strategy to pursue the minimization problem is to choose a different solution compared to the simulations of Figures \ref{fig:sim_QReq} and \ref{fig:sim_Q10}.

\section{Concluding Remarks and Future Work}\label{conc}
In this paper, we pursued optimization with primal-dual dynamics within a mean-field game. We did that by creating a network of agents with scale-free topology over which agents participating in the mean-field game communicate regarding their mean-field objectives. Results were developed on the dynamics of the mean-field game in the first, while in the second theorem the existence of a mean-field equilibrium within the mean-field game was elaborated. In addition to this, a numerical study was conducted in the form of a simulation of the mean-field game. What was found is that consensus is obtained in cases when the penalties on state deviation and control are equal, and in cases where both penalties differ significantly from each other.

Suggestions for future work can be found in the exploration of similar dynamics, however taking different network topologies into account. Also, mean-field dynamics with implemented (online) forecasting algorithms based on regression learning models will be considered \cite{boyd2011distributed, BN1}. 

\appendix
\section{Proof of Lemma 5.1}
To prove the optimal control condition as provided in \eqref{eq:optcontrol}, we firstly need to write the Hamiltonian. 
\begin{equation}
\begin{array}{lll}
    \mathcal{H}(x,\partial_x\sigma_k^T,m_k)=\\
    \inf_v \Big\{ \frac{1}{2} v^TR v +\frac{1}{2}(\rho_k-x)^TQ(\rho_k-x)\\+\partial_x\sigma_k^T(Ax+Bv+C) \Big \}.
\end{array}
\end{equation}
Then differentiating with respect to $v$ gives:
$$
Rv+B^T\partial_x\sigma_k(x,t)=0.
$$
Solving for $v$ yields the expression for $v^*_k$ from \eqref{eq:optcontrol}. Given the fact that we do not consider disturbances for now, the Hamiltonian gives us the robust mean-field dynamics:
\begin{equation}
\begin{cases}
    \partial_t\sigma_k(x,t)+H(x,\partial_x\sigma_k(x,t),\overline{m}_k) = 0 \in \mathbb{R}, \\
    \sigma_k(x,T)=\frac{1}{2}(\rho_k-x)^TS(\rho_k-x)\in \mathbb{R}, \\
    \partial_tm_k+\partial_x(m_t\partial_pH(x,\partial_x\sigma_k(x,t),m_k) = 0 \in \mathbb{R}, \\
    m_k(x,0) = m_{k0} \in \mathbb{R}.
    \label{eq:proofth1}
    \end{cases}
\end{equation}
$m_k(x,0)$ and $\sigma_k(x,T)$ are derived from the Hamilton-Jacobi-Isaacs equation. Note that the system of equations above contains two coupled partial differential equations (PDEs). Now substituting the optimal value for $v_k$, termed $v^*_k$, in the Hamiltonian, we obtain:
\begin{equation}
    \begin{split}
        \mathcal{H}&(x,\sigma_k(x,t),m_k)  = \\
         = &\frac{1}{2}\left[ x^TQx+v_k^{*^T}Rv_k^* \right]+(\partial_x\sigma_k(x,t))^T(Ax+Bv^*_k+C) \\
         = &-\frac{1}{2}(\partial_x\sigma_k(x,t))^T\left[ BR^{-1}B^T \right]\partial_x\sigma_k(x,t)\\
         +&(\partial_x\sigma_k(x,t))^T (Ax+C)+\frac{1}{2}(\rho_k-x)^TQ(\rho_k-x).
        \end{split}
\end{equation}
Using this expression for the Hamiltonian in \eqref{eq:proofth1}, we obtain the Hamilton-Jacobi-Bellman equation in \eqref{eq:systemth1}. Continuing the proof, the third equation from \eqref{eq:systemth1} can be obtained by substituting optimal control \eqref{eq:optcontrol} in \eqref{eq:proofth1}. This concludes the proof. $\quad\square$
\section{Proof of Theorem \ref{thm1}}
We start the proof by isolating the Hamilton-Jacobi-Bellman part with fixed $\rho_k$ from \eqref{eq:systemth1} and we obtain:
\begin{equation}
    \begin{array}{lll}
         \partial_t\sigma_k(x,t)-\frac{1}{2}(\partial_x\sigma_k(x,t))^T \left[BR^{-1}B^T \right] \partial_x\sigma_k(x,t)\\
         +(\partial_x\sigma_k(x,t))^T(Ax+C)\\
         +\frac{1}{2}(\rho_k-x)^TQ(\rho_k-x)=0 \\
          \in \mathbb{R}\times [0,T[,\quad\\
          \sigma_k(x,T)=\psi(\rho_k(T),x)\in\mathbb{R}. \label{eq:th2hji}
 \end{array}
\end{equation}
Now consider the following value function:
\begin{equation}
    \sigma_k(x,t)=\frac{1}{2}x^T\Phi(t)x+H(t)^Tx+\chi(t).
\end{equation}
Then \eqref{eq:th2hji} is rewritten as follows:
\begin{equation}
    \begin{array}{lll}
        \frac{1}{2}x^T\dot{\Phi}(t)x+\dot{H}(t)^Tx+\dot{\chi}(t)\\
        -\frac{1}{2}\left[\Phi(t)^Tx+H(t) \right]^T\left[BR^{-1}B^T \right]\left[\Phi(t)^Tx+H(t) \right]\\
        +\left[\Phi(t)^Tx+H(t) \right]^T(Ax+C)\\
        +\frac{1}{2}(\rho_k-x)^TQ(\rho_k-x)=0\in\mathbb{R}\times[0,T[.
 \end{array}
\end{equation}
Since the equation above is an identity in $x$ it reduces to the following set of three equations:
\begin{equation}
    \begin{cases}
    \begin{aligned}
    \dot{\Phi}(t)&+A^T\Phi(t)+\Phi^T\left[-BR^{-1}B^T\right]\Phi(t)+Q=0 \\
    &\in [0,T[,\quad \Phi(T)=S, \end{aligned}\\
    \begin{aligned}
      \dot{H}(t)&-2\Phi(t) BR^{-1}B^T H(t)+A^TH(t)\\
      &+\Phi(t)^TC-Q\rho_k(t)=0 \in [0,T[,\\
      &\quad H(T) = -S\rho_k(T),
    \end{aligned}\\
    \begin{aligned}
    \dot{\chi}&+H(t)^T\left[BR^{-1}B^T\right]H(t)+H(t)^TC\\
    &+\frac{1}{2}\rho_k(t)^TQ\rho_k(t)=0\in [0,T[,\quad \chi(T)=\frac{1}{2}\rho_k^2(T)S.
    \end{aligned}
    \label{eq:3eq}
    \end{cases}
\end{equation}
Furthermore, the optimal control is given by:
\begin{equation} \label{eq:vstark}
v_k^*=-R^{-1}B^T(\Phi(t)^Tx+H(t)). 
\end{equation}
The existence of a solution for \eqref{eq:3eq} is assured by the standard convexity-concavity assumptions, which justifies the choice for the quadratic value function \cite{tannenbaum1993h}. 

Lastly, we average the expression for optimal control and by substitution in $\tfrac{d}{dt}\overline{m}_k(t)=A\overline{m}_k+B\overline{v}^*_k+C$, we obtain the expression 
\begin{equation}
\dot{\overline{m}}_k(t)=\left[A-BR^{-1}B^T\Phi(t)\right]\overline{m}_k(t)-BR^{-1}B^TH(t)+C,
\end{equation}

which is the second equation from \eqref{eq:th2mfe}.
Now considering the stationary case with $T \rightarrow \infty$, we obtain the following set of equations:
$$\begin{aligned}
A^T\Phi+\Phi^T\left[-BR^{-1}B^T\right]\Phi+Q&=0, \\ 
-2 \Phi BR^{-1}B^T H+A^TH+\Phi^TC-Q\rho_k&=0, \\ 
H^T\left[BR^{-1}B^T\right]H+H^TC+\frac{1}{2}\rho_k^TQ\rho_k&=0. 
\end{aligned}$$
Solving for $H$ yields the following expression:
\begin{equation} \label{eq:analyticH}
H = \left[A^T-2\Phi BR^{-1}B^T \right]^{-1} \left[Q\rho_k -\Phi^TC\right].    
\end{equation}
We now substitute \eqref{eq:analyticH} in \eqref{eq:vstark} and we obtain:
\begin{equation} \label{eq:930}
\begin{split}
v^*_k&=-R^{-1}B^T\bigg( \Phi^Tx\\
&+\left[A^T-2\Phi BR^{-1}B^T \right]^{-1} \left[Q\rho_k -\Phi^TC\right] \bigg).
\end{split}
\end{equation}
Note that the last term of \eqref{eq:930} contains inverse of a sum of matrices. We may use Woodbury's matrix identity for rewriting purposes \cite{hager1989updating}. However, this results in an even longer expression. For simplification purposes, we perform the substitution 
$$
\left[A^T-2\Phi BR^{-1}B^T \right]^{-1}:=\widetilde{Z}.
$$
As a result, \eqref{eq:930} can be rewritten as:
\begin{equation}
\begin{aligned}
v^*_k&=-R^{-1}B^T\left( \Phi^Tx+\widetilde{Z} \left[ Q\rho_k -\Phi^TC\right]\right) \\
&=-R^{-1}B^T\left( \Phi^Tx-\widetilde{Z} \Phi^TC+ \widetilde{Z} Q\rho_k \right).
\end{aligned}
\end{equation}
Averaging over population $k$ we obtain that the mean state evolves according to 
$$
\bar v_k^*= -R^{-1}B^T\left( \Phi^T\overline{m}_k- \widetilde{Z} \Phi^TC+\widetilde{Z}Q
\rho_k
\right).
$$
This implies that local the interaction of neighbors consists of local averaging and local adjustment. 

\begin{equation}
\begin{array}{lll}
\dot{\overline{m}}_k(t)  = A \overline{m}_k(t) \\
-B R^{-1}B^T\left( \Phi^T\overline{m}_k- \widetilde{Z} \Phi^TC+\widetilde{Z}Q
\rho_k
\right) +C\\
=A \overline{m}_k(t) +B R^{-1}B^T \Big( \Phi^T (\rho_k - \overline{m}_k) + \widetilde{Z} \Phi^TC \\+ (-\widetilde{Z}Q - \Phi^T) \rho_k
\Big) 
+C \\
=A \overline{m}_k(t) +B R^{-1}B^T \Phi^T (\rho_k - \overline{m}_k) \\
+ B R^{-1}B^T \Big( \widetilde{Z} \Phi^TC + (-\widetilde{Z}Q - \Phi^T) \rho_k \Big) 
+C.
\end{array}
\end{equation}

Isolating the dynamics for the \textit{here and now} $c$ vector components we obtain
\begin{equation}
\begin{array}{lll}
\dot{\overline{m}}_k^c(t)  
=-\tilde Q_1 \overline{m}_k^c(t) +R_c^{-1} \Phi_c^T (\rho_c^k - \overline{m}_k^c) \\
+R_c^{-1} \Big( \widetilde{Z}_c \Phi_c^T \mu + (-\widetilde{Z}_cQ_c - \Phi_c^T) \rho_c^k \Big) 
+(\mu-\tilde f_1).
\end{array}
\end{equation}

For all aggregate states $\overline{m}^c = \left( \overline{m}_1^c,\overline{m}_2^c,\ldots,\overline{m}_p^c \right)$, we have the following consensus-type dynamics:
\begin{equation}
\dot{\overline{m}}^c = - \mathbf{\tilde Q_1} \overline{m}^c(t) -L \overline{m}^c(t) +\delta.
\end{equation}
Here, $\mathbf{\tilde Q_1}:=diag(\tilde Q_1)$ (diagonal matrix with block entry $\tilde Q_1$), $\delta:=R_c^{-1} \Big( \widetilde{Z}_c \Phi_c^T \mu + (-\widetilde{Z}_cQ_c - \Phi_c^T) \rho_c^k \Big) 
+(\mu-\tilde f_1)$ and $L$ is defined as the graph-Laplacian matrix where the $kj$th entry is the block matrix:
\begin{equation}
L_{kj} =
    \begin{cases}
    R_c^{-1} \Phi_c^T  & \quad j=k, \\
\frac{1}{|N(k)|}R_c^{-1} \Phi_c^T & \quad j \neq k, j \in N(k), \\
    0 & \quad otherwise.
    \end{cases}
\end{equation}

This concludes the proof of \textit{Theorem 5.1}. $\quad\square$                                         

\bibliographystyle{plain}        
\bibliography{rolingbauso}           

\begin{thebibliography}{10}

\bibitem{achdou2012mean}
Yves Achdou, Fabio Camilli, and Italo Capuzzo-Dolcetta.
\newblock Mean field games: numerical methods for the planning problem.
\newblock {\em SIAM Journal on Control and Optimization}, 50(1):77--109, 2012.

\bibitem{bauso2017consensus}
Dario Bauso.
\newblock Consensus via multi-population robust mean-field games.
\newblock {\em Systems \& Control Letters}, 107:76--83, 2017.

\bibitem{bauso2017dynamic}
Dario Bauso.
\newblock Dynamic demand and mean-field games.
\newblock {\em IEEE Transactions on Automatic Control}, 62(12):6310--6323,
  2017.

\bibitem{bauso2016crowd}
Dario Bauso, Thulasi Mylvaganam, and Alessandro Astolfi.
\newblock Crowd-averse robust mean-field games: approximation via state space
  extension.
\newblock {\em IEEE Transactions on Automatic Control}, 61(7):1882--1894, 2016.

\bibitem{BN1}
Dario Bauso and Toru Namarikawa.
\newblock {\em 20. Electric vehicles and Mean-field, from Advanced Data
  Analytics for Power Systems edited by Ali Tajer; Samir M. Perlaza; H. Vincent
  Poor}.
\newblock Cambridge University Press, in print.

\bibitem{bertsimas2000traffic}
Dimitris Bertsimas and Sarah~Stock Patterson.
\newblock The traffic flow management rerouting problem in air traffic control:
  A dynamic network flow approach.
\newblock {\em Transportation Science}, 34(3):239--255, 2000.

\bibitem{boyd2011distributed}
Stephen Boyd, Neal Parikh, and Eric Chu.
\newblock {\em Distributed optimization and statistical learning via the
  alternating direction method of multipliers}.
\newblock Now Publishers Inc, 2011.

\bibitem{chen2011convergence}
Junting Chen and Vincent~KN Lau.
\newblock Convergence analysis of saddle point problems in time varying
  wireless systems—control theoretical approach.
\newblock {\em IEEE Transactions on Signal Processing}, 60(1):443--452, 2011.

\bibitem{chen2019exponential}
Xin Chen and Na~Li.
\newblock Exponential stability of primal-dual gradient dynamics with
  non-strong convexity.
\newblock {\em arXiv preprint arXiv:1905.00298}, 2019.

\bibitem{cherukuri2017saddle}
Ashish Cherukuri, Bahman Gharesifard, and Jorge Cortes.
\newblock Saddle-point dynamics: conditions for asymptotic stability of saddle
  points.
\newblock {\em SIAM Journal on Control and Optimization}, 55(1):486--511, 2017.

\bibitem{cherukuri2016asymptotic}
Ashish Cherukuri, Enrique Mallada, and Jorge Cort{\'e}s.
\newblock Asymptotic convergence of constrained primal--dual dynamics.
\newblock {\em Systems \& Control Letters}, 87:10--15, 2016.

\bibitem{cherukuri2017role}
Ashish Cherukuri, Enrique Mallada, Steven Low, and Jorge Cort\'es.
\newblock The role of convexity in saddle-point dynamics: Lyapunov function and
  robustness.
\newblock {\em IEEE Transactions on Automatic Control}, 63(8):2449--2464, 2017.

\bibitem{feijer2010stability}
Diego Feijer and Fernando Paganini.
\newblock Stability of primal--dual gradient dynamics and applications to
  network optimization.
\newblock {\em Automatica}, 46(12):1974--1981, 2010.

\bibitem{ferragut2014network}
Andr{\'e}s Ferragut and Fernando Paganini.
\newblock Network resource allocation for users with multiple connections:
  fairness and stability.
\newblock {\em IEEE/ACM Transactions on Networking (TON)}, 22(2):349--362,
  2014.

\bibitem{hager1989updating}
William~W Hager.
\newblock Updating the inverse of a matrix.
\newblock {\em SIAM review}, 31(2):221--239, 1989.

\bibitem{halldorsson2010sustainable}
{\'A}rni Halld{\'o}rsson and Gy{\"o}ngyi Kov{\'a}cs.
\newblock The sustainable agenda and energy efficiency: Logistics solutions and
  supply chains in times of climate change.
\newblock {\em International Journal of Physical Distribution \& Logistics
  Management}, 40(1/2):5--13, 2010.

\bibitem{jonsson2010primal}
Ulf~T J{\"o}nsson.
\newblock Primal and dual criteria for robust stability and their application
  to systems interconnected over a bipartite graph.
\newblock In {\em Proceedings of the 2010 American Control Conference}, pages
  5458--5464. IEEE, 2010.

\bibitem{li2017distributed}
Chaojie Li, Xinghuo Yu, Tingwen Huang, and Xing He.
\newblock Distributed optimal consensus over resource allocation network and
  its application to dynamical economic dispatch.
\newblock {\em IEEE Transactions on Neural Networks and Learning Systems},
  29(6):2407--2418, 2017.

\bibitem{li2017dynamic}
Mushu Li, Peter He, and Lian Zhao.
\newblock Dynamic load balancing applying water-filling approach in smart grid
  systems.
\newblock {\em IEEE Internet of Things Journal}, 4(1):247--257, 2017.

\bibitem{liang2019exponential}
Shu Liang, Le~Yi Wang, and George Yin.
\newblock Exponential convergence of distributed primal--dual convex
  optimization algorithm without strong convexity.
\newblock {\em Automatica}, 105:298--306, 2019.

\bibitem{montreuil2011toward}
Benoit Montreuil.
\newblock Toward a physical internet: meeting the global logistics
  sustainability grand challenge.
\newblock {\em Logistics Research}, 3(2-3):71--87, 2011.

\bibitem{montreuil2010towards}
Benoit Montreuil, Russell~D Meller, and Eric Ballot.
\newblock Towards a physical internet: the impact on logistics facilities and
  material handling systems design and innovation.
\newblock In {\em Proceedings of the 11th edition of the IMHRC in Milwaukee,
  Wisconsin, USA - 2010}, 2010.

\bibitem{nguyen2018contraction}
Hung~D Nguyen, Thanh~Long Vu, Konstantin Turitsyn, and Jean-Jacques Slotine.
\newblock Contraction and robustness of continuous time primal-dual dynamics.
\newblock {\em IEEE control systems letters}, 2(4):755--760, 2018.

\bibitem{notarstefano2019distributed}
Giuseppe Notarstefano, Ivano Notarnicola, and Andrea Camisa.
\newblock Distributed optimization for smart cyber-physical networks.
\newblock {\em Foundations and Trends in Systems and Control}, 7(3):253--383,
  2019.

\bibitem{piecyk2015green}
Maja Piecyk, Michael Browne, Anthony Whiteing, and Alan McKinnon.
\newblock {\em Green logistics: Improving the environmental sustainability of
  logistics}.
\newblock Kogan Page Publishers, 2015.

\bibitem{stolyar2006greedy}
Alexander~L Stolyar.
\newblock Greedy primal-dual algorithm for dynamic resource allocation in
  complex networks.
\newblock {\em Queueing Systems}, 54(3):203--220, 2006.

\bibitem{tannenbaum1993h}
Allen Tannenbaum.
\newblock $\mathbb{H}^\infty$-optimal control and related minimax design
  problems (tamer basar and pierre bernhard).
\newblock {\em SIAM Review}, 35(3):538--540, 1993.

\bibitem{yang2018distributed}
Chungang Yang, Haoxiang Dai, Jiandong Li, Yue Zhang, and Zhu Han.
\newblock Distributed interference-aware power control in ultra-dense small
  cell networks: a robust mean field game.
\newblock {\em IEEE Access}, 6:12608--12619, 2018.

\bibitem{yang2016distributed}
Chungang Yang, Jiandong Li, Prabodini Semasinghe, Ekram Hossain, Samir~M
  Perlaza, and Zhu Han.
\newblock Distributed interference and energy-aware power control for
  ultra-dense d2d networks: A mean field game.
\newblock {\em IEEE Transactions on Wireless Communications}, 16(2):1205--1217,
  2016.

\bibitem{zhao2012swing}
Changhong Zhao, Ufuk Topcu, and Steven Low.
\newblock Swing dynamics as primal-dual algorithm for optimal load control.
\newblock In {\em 2012 IEEE Third International Conference on Smart Grid
  Communications (SmartGridComm)}, pages 570--575. IEEE, 2012.

\end{thebibliography}

\end{document}